\documentclass[amstex,12pt,russian,amssymb]{article}
\usepackage{mathtext}
\usepackage[cp1251]{inputenc}
\usepackage[russian]{babel}
\usepackage{amsmath}
\usepackage{amssymb}
\usepackage{amsxtra}
\usepackage{latexsym}
\usepackage{ifthen}

\usepackage{color}

\newcommand {\beq}{\begin{equation}}
\newcommand {\eeq}{\end{equation}}
\newcommand{\reff}[1]{(\ref{#1})}

\newtheorem{theorem}{Theorem}[section]
\newtheorem{lemma}[theorem]{Lemma}

\newtheorem{definition}[theorem]{Definition}
\newtheorem{remark}[theorem]{Remark}

\begin{document}

\title{A local large deviation principle\\ for inhomogeneous birth-death processes}

\author{N.D. Vvedenskaya~$^{1}$, A.V. Logachov~$^{2,3,4}$, Y.M. Suhov~$^{1,5}$,
 A.A. Yambartsev~$^{6}$ }

\maketitle

 {\footnotesize

 \noindent $^1$ Dobrushin Laboratory, Institute for Information Transmission Problems, RAS, 19 Bol'shoi
Karetnyi per, Moscow 127051, RF\\
E-mail:\ \ ndv@iitp.ru\footnote{The work by N.D. Vvedenskaya had been conducted at the IITP RAS
under the support of the RSF (Project  No14-50-00 150).}

\noindent $^2$ Laboratory of Applied Mathematics, Novosibirsk State University, 2 Pirogov ul,
Novosibirsk 630090, RF\\

\noindent $^3$ Laboratory of Probability Theory and Mathematical Statistics, Sobolev Institute of Mathematrics,
Siberiian Branch of the RAS, 4 Koptyug ul, Novosibirsk 630090, RF

\noindent $^4$ Statistics Division, Novosibirsk State University of Economics and Management,
56 Kamenskaya ul, Novosibirsk 630099, RF\\
E-mail:\ \ omboldovskaya@mail.ru\footnote{The work by A.V. Logachov was supported by the RSF Grant
(Project No18-11-00129). }

\noindent $^5$  Math Dept, Penn State University, University Park, State College, PA 16802, USA.\\
E-mail:\ \ yms@statslab.cam.ac.uk

\noindent $^6$ Department of Statistics, Institute of Mathematics and Statistics, University of S\~ao Paulo,
1010 Rua do Mat\~ao, CEP 05508--090, S\~ao Paulo SP, Brazil \\
E-mail:\ \ yambar@gmail.com}


\vspace{1cm}

\begin{abstract}
The paper considers a continuous-time birth-death process where the jump rate has an asymptotically
polynomial dependence on the process position. We obtain a rough exponential asymptotics
for the probability of excursions of a re-scaled process contained within a neighborhood of a given
continuous non-negative function.
\end{abstract}

\section{Introduction} \def\lam{\lambda}

In the modern literature on the large deviation principle, one considers various conditions for random
processes guaranteeing a rough exponential asymptotics for probabilities of rare events. See,
for example, \cite{DZ} -- \cite{Var}.
In this paper we deal with birth-and-death Markov processes that are inhomogeneous in the state space:
the rates of jumps are polynomially dependent on the position of the process. For these processes we obtain
exponential asymptotics for the probabilities the normalized process to be in a neighborhood of a continuous
function. Moreover, we provide this asymptotics both for ergodic processes and for transient (even exploding)
processes.

The study of birth-and-death processes is of a certain mathematical interest and, moreover, is important
for a number of applications. As examples, we can cite the information theory (encoding and storage of
information, see \cite{SS1}), biology and chemistry (models of growth and extinction in systems with multiple
components,  see \cite{SS2},  \cite{MSSZ}), and economics (models of competitive production and pricing,
\cite{MPY}, \cite{VSB}).

Consider a continuous-time Markov process $\xi (t)$, $t\geq 0$, with state space $\mathbb{Z}^+\cup\{\infty\}$,
where $\mathbb{Z}^+=\{0\}\cup\mathbb{N}$. Let us assume that the process starts at 0.

The evolution of the process $\xi (\,\cdot\,)$ is described as follows. For a given $t\ge 0$, let
$\xi (t)=x\in\mathbb{Z}^+$. The state of the random process does not change during the random time
$\tau_x$ with exponential distribution  with parameter $h(x)>0$. At the moment $t+\tau_x$ the process
jumps to the states $x\pm 1$ with probabilities
\beq
{\bf P}(\xi (t+\tau_x)=x+1)=\frac{\lambda(x)}{h(x)},\ \ \ {\bf P}(
  \xi(t+\tau_x)=x-1)=\frac{\mu(x)}{h(x)}, \label{2}
\eeq
correspondingly, where $\lambda(x)+\mu(x)=h(x)$, $\lambda(x)> 0$ when $x\in \mathbb{Z}^+$, and
$\mu(x)> 0$ for $x\in \mathbb{N}$.

Suppose that for $x=0$ the rates $\mu (x)=0$, $\lambda(x)=\lambda_0>0$ (i.e. the process cannot take
negative values), and the following asymptotics hold true
\beq
\lim_{x\to
\infty}\frac{\lambda(x)}{P_lx^l}= \lim_{x\to
\infty}\frac{\mu(x)}{Q_mx^m}=1, \label{1}
\eeq
where $P_l$ and $Q_m$ are positive constants, and $l\geq 0$,
$m\geq 0$, $l\vee m>0$.

When $l\leq 1$, the  existence of a Markov process with the above properties
is established in the standard way, see for example \cite{Fell}, Ch. 17, \S 4, 5, and also  \cite{KS},
Ch. 2, $\S$ 5, Theorem 2.5.5, \cite{KT}, Ch. 6, 7. When $l>1$ the process $\xi(\cdot)$, generally
speaking, can go to infinity ("explode") during a random time, finite with probability one. There are
two approaches to construct such processes. (1) One can stop the process at the random time point
(the time of explosion); viz., see \cite{Korol}, Ch. 15, \S 4,
\cite{Dyn},  vol. {\bf 1}, Ch 6, P. 365; vol. {\bf 2},  P. 274.  \  (2)\ One can extend the phase space
$\mathbb{Z}^+$ by adding an
absorbing state (denoted by $\infty$); see, e.g., \cite{Ito}, Ch. 4, \S 48, \cite{Fell}, Ch. 17, $\S$ 10.
In this paper we use the second approach.

The above class of random processes has been given the name birth-and-death processes; see,
for example,  \cite{Fell},  \cite{Ito}.

There exist conditions on $l$ and $m$ which are sufficient for explosion and
non-explosion. For example, when $l>1$ and $m<l$, the process $\xi (\,\cdot\,)$ explodes, while if
$m>l$ it does not.  As references, cf. original papers
\cite{KMcG}, \cite{LR} and references within. See also \cite{No}, Ch. 23, $\S$ 7,
\cite{KS}, Ch. 2, $\S$ 5, and \cite{St}, Ch.5, $\S$ 3 (the last reference includes results for general Markov
chains, not only for birth-and-death processes).

We are interested in the local large deviation principle (LLDP) for the family of scaled processes
\begin{equation}\label{xit}
\xi_T(t)=\frac{\xi(tT)}{T}, \ 0\leq t\leq 1,
\end{equation}
where $T>0$ is a parameter  (see, \cite{Bor-Mog2}, \cite{Bor-Mog1}).
In a sense, the formulation and analysis of the LLDP should precede the study of
other forms of the large deviation principle.

The validity of our results does not depend on whether or not the process $\xi(\cdot )$ explodes within
a finite time. We focus on the asymptotics of the probability of the event that the trajectories of process
$\xi_T(\,\cdot\, )$ to stay in a neighborhood of a continuous positive function given on the interval $[0,1]$.
It means that we are working on the set of trajectories which do not tend to infinity in the time interval
$[0,T]$. The considered probabilities are positive even if the process $\xi(\cdot)$ explodes (see equation
(\ref{8a}) below).

Let $\mathbb{D}[0,1]$ denote the space of right-continuous functions with left-limit at eact $t\in [0,1]$.
For any  $f,g \in \mathbb{D}[0,1]$, set
$$
\rho(f,g)=\sup\limits_{t\in[0,1]}|f(t)-g(t)|.$$
\def\diy{\displaystyle}

\begin{definition}\label{def1}
The family of random processes $\xi_T(\,\cdot\,)$ satisfies the
LLDP on the set $G\subseteq \mathbb{D}[0,1]$ with a rate functional
$I = I(f)\,:\, \mathbb{D}[0,1] \rightarrow [0,\infty)$ and a normalising function
$T\in (0,\infty )\mapsto \psi(T)>0$ with $\lim\limits_{T\rightarrow\infty}\psi(T) = \infty$
if, for any function $f \in G$, the following equality holds true:
\begin{equation}\label{opr}\begin{array}{l}
\lim\limits_{\varepsilon\rightarrow 0}\limsup\limits_{T\rightarrow \infty}\diy\frac{1}{\psi(T)}
\ln\mathbf{P}(\xi_T(\,\cdot\,)\in U_\varepsilon(f))\\
\qquad =\lim\limits_{\varepsilon\rightarrow 0}\liminf\limits_{T\rightarrow \infty}\diy\frac{1}{\psi(T)}
\ln\mathbf{P}(\xi_T(\,\cdot\,)\in U_\varepsilon(f))=-I(f).\end{array}\end{equation}
Here
\begin{equation}\label{oprU}
U_\varepsilon(f)=\{g\in \mathbb{D}[0,1]: \ \rho(f,g)<\varepsilon\}.\end{equation}
\end{definition}

In the framework of Definition~\ref{def1} there are various cases to consider.
We separate three cases: 1) $l>m$, 2) $l<m$ and $l=m$.

Note that the case $m=1$, $l=0$, follows from \cite{MPY}
(where a two dimensional Markov process is treated). A
similar result is obtained in \cite{Log1}
for solutions of stochastic differential Ito's equations.
The classical case $l=m=0$, $\varphi(T)=T$ follows, for example, from \cite{Bor-Mog}.

In this paper we use the approach developed in \cite{MPY}. We would like to note that the
large deviation principle for the sequence of processes $\xi_T(\,\cdot\,)$ in space
$\mathbb{D}[0,1]$ with Skorohod metric cannot be obtained even for non-exploding processes:
one can show that the corresponding family of
measures is not exponentially dense, except for the case $l=m=0$.

The paper is organized as follows: in $\S$ 2 we introduce our definitions and the system
of notation, as well as the main result (Theorem~\ref{th2.1})
and key lemmas. In $\S$ 3 we prove Theorem~\ref{th2.1} and key lemmas.
In $\S$ 4 (the Appendix) some auxiliary technical assertions are established.


\section{Main results, definitions}

Let $F$ denote the set of functions $f(t) \in \mathbb{C}[0,1]$ such that  $f(0)=0$ and
$f(t)>0$ as $0<t\leq 1$.

\begin{theorem} \label{th2.1}
Let the conditions \eqref{2} and \eqref{1} be fulfilled. Than the random processes sequence
 $\xi_T(\cdot)$ on  $F$ fulfills  the following LLDP:

{\rm{a)}} If  $l>m$ then the normalizing function  $\psi(T)=T^{l+1}$,  and the rate functional has the form
$$
I(f)=P_l\int_0^1f^l(t)dt, \ \ \ f\in F.
$$

{\rm{b)}} If $l=m$ and  $P_l\neq Q_m$ then $\psi(T)=T^{l+1}$ and
$$
I(f)=(\sqrt{P_l}-\sqrt{Q_m})^2\int_0^1f^l(t)dt, \ \ \ f\in F.
$$

{\rm{c)}} If $l<m$ then  $\psi(T)=T^{m+1}$ and
$$
I(f)=Q_m\int_0^1f^m(t)dt, \ \ \ f\in F.
$$
\end{theorem}

The case where $l=m$ and $P_l= Q_m$ needs a different normalization; we do not discuss  it
in this paper.


\bigskip

Consider a space- and time-homogenous Markov process
$\zeta(t)$, $t\in[0,T]$, on the phase space $\mathbb{Z}$, where the jump rate is
equal to $1$, and the jump size is  $\pm 1$, occurring with probability $1/2$.

Denote by $X_T$ the set of all right-continuous step-functions with a finite number of
$\pm 1$-jumps on  $[0,T]$.

\begin{lemma} \label{l2.0}

For any given  $T$, the distribution  ${\mathbf P}^{(\xi)}_T(\,\cdot\,\cap X_T)$
of process  $\xi(\cdot)$ on $X_T$ is absolutely continuous with respect to the distribution
${\mathbf P}^{(\zeta )}_T$
of process $\zeta(\cdot)$ on $X_T$. The corresponding density (the Radon-Nikodym  derivative)
${\mathfrak p}_T(u)=\diy\frac{d{\mathbf P}^{(\xi )}_T}{d{\mathbf P}^{(\zeta )}_T}(u)$, $u\in X_T$,
has the form
\beq
{\mathfrak p}_T(u)=\left\{ \begin{array}{ll}2^{N_T(u)}\left(\prod\limits_{i=1}^{N_T(u)}
e^{-(h(u(t_{i-1}))-1)\tau_{i}}\nu(u(t_{i-1}),u(t_i))\right)\\
\qquad \times e^{-(h(u(t_{N_T(u)})-1))(T-t_{N_T(u)})
},\;\;\qquad\quad\quad \;\; \mbox{ if } \;N_T(u)\geq 1,\\
e^{-(h(0)-1)T}, \qquad\qquad\qquad\qquad\qquad\qquad\qquad  \mbox{ if } \;  N_T(u)= 0.\\
\end{array} \right. \label{8a}\eeq
Here it is supposed that  the function $u(\cdot)$ on  $[0,T]$ has exactly  $N_T(u)$ jumps
at the time points   \
$t_1,t_2,...,t_{N_T(u)}$ \ where  \ $0=t_0<t_1<...<t_{N_T(u)}\leq T$, \
$\tau_i=t_i-t_{i-1}$. Moreover,
$$
\nu (u(t_{i-1}),u (t_i))= \left\{ \begin{array}{ll}\lambda(u (t_{i-1})), &
\mbox{ if }\; u (t_i)-u (t_{i-1})=1;\\
\mu(u (t_{i-1})), & \mbox{ if }\;  u (t_i)-u (t_{i-1})= -1. \end{array} \right.
$$

\end{lemma}

Observe that the probability ${\mathbf P}(\xi(\cdot)\in X_T)$ in Lemma  2.2
is allowed to be less then 1. (Clearly, this probability is positive.) The same density ${\mathfrak p}_T$
was used in  \cite{MPY}.

Let us denote by $N_T(\zeta)$  the random number of jumps in process $\zeta(\cdot)$ on interval
$[0,T]$.

The assertion of Lemma \ref{l2.0} is equivalent to the fact that for any measurable set  $G\subseteq X_T$
\beq \label{3}
\mathbf{P}(\xi(\cdot)\in G)=e^T \mathbf{E}(e^{-A_T(\zeta)}e^{B_T(\zeta)+N_T(\zeta)\ln2};\zeta(\cdot)\in G).
\eeq
We set
\beq  \label{4}\begin{array}{l}
A_T(\zeta)=\int_0^T h(\zeta(t))dt\\
\quad =\left\{ \begin{array}{ll}\sum\limits_{i=1}^{N_T(\zeta)}h(\zeta(t_{i-1}))\tau_{i}+
h(\zeta(t_{N_T(\zeta)}))(T-t_{N_T(\zeta)}),&
\mbox{if }\; N_T(\zeta)\geq1,\\
h(0)T, & \mbox{ if }\;  N_T(\zeta)=0; \end{array} \right.\\
B_T(\zeta)=\left\{ \begin{array}{ll}\sum\limits_{i=1}^{N_T(\zeta)}\ln(\nu(\zeta(t_{i-1}),\zeta(t_i))), &
\mbox{if }\; N_T(\zeta)\geq1,\\
0, & \mbox{ if }\;  N_T(\zeta)=0. \end{array} \right.\end{array}\eeq
The expressions in (\ref{3}) specify, in our context, the statement of the Radon-Nikodym theorem
(see, e.g.,  \cite{D-Sch}, Theorem 2, sec. III, ch. 10). The expressions (\ref{3}) are used for
analysing the asymptotical behaviour of the logarithm of probability
$\mathbf{P}(\xi_T(\cdot)\in U_\varepsilon(f))$, $f\in F$.

Theorem \ref{th2.1} indicates that for $l\neq m$  the main contribution into  the
asymptotics is brought by
 $A_T(\zeta)$, whereas in the case $l=m$ the asymptotics involves both
 $A_T(\zeta)$ and $B_T(\zeta)$.

Consider the family of scaled processes
$$\zeta_T(t)=\frac{\zeta(tT)}{T}, \ \ t\in[0,1].$$
Let $k_+$ and  $k_-$ denote the number of positive and negative jumps in  $\zeta_T(\cdot)$  and
set $L=k_+-k_-$.

For $\zeta_T(\cdot)\in U_\varepsilon(f)$ we have the inequality
\beq
f(1)-\varepsilon \leq \zeta_T(1) \leq f(1)+\varepsilon.  \label{8}
\eeq
The jumps in $\zeta_T(\cdot)$  are  $\pm 1/T$, therefore
 (\ref{8}) yields the inequalities
\beq
(f(1)-\varepsilon)T \leq L \leq (f(1)+\varepsilon)T.  \label{9}
\eeq
With these definitions and observations we can write:
\beq
k_+ + k_- =N_T(\zeta), \ \ \  k_+=\frac{N_T(\zeta)+L}{2}, \ \ \ k_-=\frac{N_T(\zeta)-L}{2}. \label{10}
\eeq
For brevity, we  write  below $\xi_T, \zeta_T$ and $A_T, B_T$ instead of  $\xi_T(\cdot), \zeta_T(\cdot)$
and $A_T(\zeta), B_T(\zeta)$. Also set: $v=\max (l,m)$.

\smallskip

\begin{lemma}\label{l2.3} Let $f\in F$. In case $l\neq m$ we have
$$\lim\limits_{\varepsilon\rightarrow 0} \limsup\limits_{T\rightarrow\infty} \frac{1}{T^{v+1}}
\ln\mathbf{E}(e^{B_T+N_T(\zeta)\ln2};\zeta_T\in U_\varepsilon(f))\leq0,$$
whereas in case $l=m$
$$\lim\limits_{\varepsilon\rightarrow 0} \limsup\limits_{T\rightarrow\infty} \frac{1}{T^{l+1}}
\ln\mathbf{E}(e^{B_T+N_T(\zeta)\ln2};\zeta_T\in U_\varepsilon(f))\leq
2\sqrt{P_lQ_m}\int_0^1f^l(s)ds,$$

\end{lemma}

\begin{lemma} \label{l2.4} For $f\in F$, in case $l\neq m$
$$\lim\limits_{\varepsilon\rightarrow 0} \liminf\limits_{T\rightarrow\infty} \frac{1}{T^{v+1}}
\ln\mathbf{E}(e^{B_T+N_T(\zeta)\ln2};\zeta_T\in U_\varepsilon(f))\geq0,$$
and in case $l=m$
$$\lim\limits_{\varepsilon\rightarrow 0} \liminf\limits_{T\rightarrow\infty} \frac{1}{T^{l+1}}
\ln\mathbf{E}(e^{B_T+N_T(\zeta)\ln2};\zeta_T\in U_\varepsilon(f))\geq 2\sqrt{P_lQ_m}\int_0^1f^l(s)ds .$$
\end{lemma}

\section{Proofs of Theorem \ref{th2.1} and Lemmas \ref{l2.0}--\ref{l2.4}}

\smallskip
In what follows, $\blacktriangle$ marks the end of a proof.

{\sc Proof of Theorem \ref{th2.1}} We are going to get the LLDP for functions  $f\in F$. First
let us estimate the quantity $A_T$. Fix a value $\varepsilon >0$ until a further notice.

From equation \eqref{4} it follows that
$$
A_T:=\int_0^T h(\zeta(t))dt=T\int_0^1 h(T\zeta_T(s))ds.
$$
If $\zeta_T\in U_\varepsilon(f)$ then
\beq
(f(s)-\varepsilon)\leq \zeta_T(s)\leq (f(s)+\varepsilon). \label{11}
\eeq

Let $\delta\in(0,1)$ be also fixed for the time being and denote
$m_\delta:=\min\limits_{t\in[\delta,1]}f(t)$. Here $m_\delta>0$ for $f\in F$. Therefore,
$k_0=m_\delta-\varepsilon>0$ when $\varepsilon$ is sufficiently small.

Let us estimate  $A_T$ on the set of trajectories where  inequality  (\ref{11}) is valid.
From (\ref{11}) it follows that $T\zeta_T(s)\geq k_0 T$ for $s\in[\delta,1]$. Therefore,
by virtue of condition (\ref{1}), for any $\gamma_0\in (0,1)$ and  $s\in[\delta,1]$,
for $T$ large enough we have the inequalities
\beq\label{12}
1-\gamma_0\leq \frac{h(T\zeta_T(s))}{P_l(T\zeta_T(s))^l}
\leq 1+\gamma_0\ \ \mbox{ in case }\ l>m,\eeq
\beq\label{13}
1-\gamma_0\leq\frac{h(T\zeta_T(s))}{(P_l+Q_m)(T\zeta_T(s))^l}
\leq 1+\gamma_0 \ \   \mbox{in case }\   l=m,\eeq
and
\beq
1-\gamma_0\leq\frac{h(T\zeta_T(s))}{Q_m(T\zeta_T(s))^m}
\leq 1+\gamma_0\ \  \mbox{ in case }\  \ l<m. \label{14}\eeq

Consider the case $l>m$. Owing to (\ref{11}) and (\ref{12}), for $T$ sufficiently large,
we get
\beq\label{15}\begin{array}{l}
T\int_\delta^1(1-\gamma_0)P_l(T(f(s)-\varepsilon))^lds
\ \leq \ A_T\\
\qquad\qquad\leq T\int_0^\delta h(T\zeta_T(s))ds+T\int_\delta^1(1+\gamma_0)P_l(T(f(s)+\varepsilon))^l
ds.\end{array}\eeq

Set $M:=\max(\max\limits_{t\in[0,1]}f(t),1)$. By using (\ref{11}), for $T$ large enough we have that
$$h(T\zeta_T(s))\leq(1+\gamma_0)P_l(T(M+\varepsilon))^l.$$
Consequently, from (\ref{15}) we obtain the inequality
\beq\begin{array}{l}
T^{l+1}P_l\int_\delta^1(1-\gamma_0)(f(s)-\varepsilon)^lds
\ \leq \ A_T \\
\qquad\leq T^{l+1}P_l\delta(1+\gamma_0)(M+\varepsilon)^l+T^{l+1}P_l\int_\delta^1(1+\gamma_0)(f(s)+\varepsilon)^lds.\end{array} \label{16}\eeq
By using the bound (\ref{16}) and equation (\ref{3}), we get the following:
\beq\label{17}\begin{array}{l}
\exp\,[-T^{l+1}P_l\int_\delta^1(1-\gamma_0)(f(s)-\varepsilon)^lds]
e^T \mathbf{E}(e^{B_T+N_T(\zeta)\ln2};\zeta_T\in U_\varepsilon(f))\\
\qquad
\geq\mathbf{P}(\xi_T(\cdot)\in U_\varepsilon(f))=e^T \mathbf{E}(e^{-A_T}e^{B_T+N_T(\zeta)\ln2};\zeta_T\in U_\varepsilon(f))\\
\qquad \geq\exp\,[-T^{l+1}P_l\delta(1+\gamma_0)(M+\varepsilon)^l-T^{l+1}P_l\int_\delta^1(1+\gamma_0)(f(s)+\varepsilon)^lds]\\
\quad\qquad\times e^T \mathbf{E}(e^{B_T+N_T(\zeta)\ln2};\zeta_T\in U_\varepsilon(f)).\end{array}\eeq
Further, by virtue of (\ref{17})
\beq\begin{array}{l}
-P_l\int_\delta^1(1-\gamma_0)(f(s)-\varepsilon)^lds +
\limsup\limits_{T\rightarrow\infty}{\diy\frac{1}{T^{l+1}}}\ln\mathbf{E}(e^{B_T+N_T(\zeta)\ln2};\zeta_T\in U_\varepsilon(f))\\
\qquad\geq\limsup\limits_{T\rightarrow\infty}{\diy\frac{1}{T^{l+1}}}\ln\mathbf{P}(\xi_T\in U_\varepsilon(f))\geq\liminf\limits_{T\rightarrow\infty}{\diy\frac{1}{T^{l+1}}}\ln\mathbf{P}(\xi_T\in U_\varepsilon(f))\\
\qquad\geq -P_l\delta(1+\gamma_0)(M+\varepsilon)^l-P_l\int_\delta^1(1+\gamma_0)(f(s)+\varepsilon)^lds\\
\qquad\quad+\liminf\limits_{T\rightarrow\infty}{\diy\frac{1}{T^{l+1}}}
\ln\mathbf{E}(e^{B_T+N_T(\zeta)\ln2};\zeta_T\in U_\varepsilon(f)).
\end{array} \label{18}
\eeq
Next, from (\ref{18}) it follows that
\beq
\begin{array}{l}
-P_l\int_\delta^1(1-\gamma_0)f^l(s)ds+
\lim\limits_{\varepsilon\rightarrow 0}\limsup\limits_{T\rightarrow\infty}{\diy\frac{1}{T^{l+1}}}\ln\mathbf{E}(e^{B_T+N_T(\zeta)\ln2};\zeta_T\in U_\varepsilon(f))\\
\quad\geq\lim\limits_{\varepsilon\rightarrow 0}\limsup\limits_{T\rightarrow\infty}{\diy\frac{1}{T^{l+1}}}\ln\mathbf{P}(\xi_T(\cdot)\in U_\varepsilon(f))\geq\lim\limits_{\varepsilon\rightarrow 0}\liminf\limits_{T\rightarrow\infty}
{\diy\frac{1}{T^{l+1}}}\ln\mathbf{P}(\xi_T(\cdot)\in U_\varepsilon(f))\\
\quad\geq-P_l\delta(1+\gamma_0)M^l-P_l\int_0^1(1+\gamma_0)f^l(s)ds\\
\qquad\qquad\qquad\qquad\qquad\qquad\quad
+\lim\limits_{\varepsilon\rightarrow 0}\liminf\limits_{T\rightarrow\infty}{\diy\frac{1}{T^{l+1}}}
\ln\mathbf{E}(e^{B_T+N_T(\zeta)\ln2};\zeta_T\in U_\varepsilon(f)).
\end{array} \label{19}
\eeq

Note that the inequatity  (\ref{19}) is valid for all $\gamma_0,\delta>0$. Letting
 $\gamma_0\rightarrow 0$ and $\delta\rightarrow 0$ we get that
\beq\begin{array}{l} -P_l\int_0^1f^l(s)ds+
\lim\limits_{\varepsilon\rightarrow 0}\limsup\limits_{T\rightarrow\infty}{\diy\frac{1}{T^{l+1}}}\ln\mathbf{E}(e^{B_T+N_T(\zeta)\ln2};\zeta_T\in U_\varepsilon(f))\\
\quad\geq\lim\limits_{\varepsilon\rightarrow 0}\limsup\limits_{T\rightarrow\infty}\frac{1}{T^{l+1}}\ln\mathbf{P}(\xi_T(\cdot)\in U_\varepsilon(f))\geq\lim\limits_{\varepsilon\rightarrow 0}\liminf\limits_{T\rightarrow\infty}
{\diy\frac{1}{T^{l+1}}}\ln\mathbf{P}(\xi_T(\cdot)\in U_\varepsilon(f))\\
\quad\geq-P_l\int_0^1f^l(s)ds+\lim\limits_{\varepsilon\rightarrow 0}\liminf\limits_{T\rightarrow\infty}
{\diy\frac{1}{T^{l+1}}}
\ln\mathbf{E}(e^{B_T+N_T(\zeta)\ln2};\zeta_T\in U_\varepsilon(f)).
\end{array} \label{20}
\eeq

In a similar way, by using (\ref{13}) and  (\ref{14}) we obtain inequalities
for the case $ l = m $:
\beq \begin{array}{l}
-(P_l+Q_m)\int_0^1f^l(s)ds+
\lim\limits_{\varepsilon\rightarrow 0}\limsup\limits_{T\rightarrow\infty}{\diy\frac{1}{T^{l+1}}}\ln\mathbf{E}(e^{B_T+N_T(\zeta)\ln2};\zeta_T\in U_\varepsilon(f))\\
\quad\geq\lim\limits_{\varepsilon\rightarrow 0}\limsup\limits_{T\rightarrow\infty}{\diy\frac{1}{T^{l+1}}}\ln\mathbf{P}(\xi_T\in U_\varepsilon(f))\geq\lim\limits_{\varepsilon\rightarrow 0}\liminf\limits_{T\rightarrow\infty}
{\diy\frac{1}{T^{l+1}}}\ln\mathbf{P}(\xi_T\in U_\varepsilon(f))\\
\quad\geq-(P_l+Q_m)\int_0^1f^l(s)ds+\lim\limits_{\varepsilon\rightarrow 0}\liminf\limits_{T\rightarrow\infty}
{\diy\frac{1}{T^{l+1}}}
\ln\mathbf{E}(e^{B_T+N_T(\zeta)\ln2};\zeta_T\in U_\varepsilon(f)),\end{array} \label{21}\eeq
and for the case $l<m$:
\beq\begin{array}{l}
-Q_m\int_0^1f^m(s)ds+
\lim\limits_{\varepsilon\rightarrow 0}\limsup\limits_{T\rightarrow\infty}{\diy\frac{1}{T^{m+1}}}\ln\mathbf{E}(e^{B_T+N_T(\zeta)\ln2};\zeta_T\in U_\varepsilon(f))\\
\quad\geq\lim\limits_{\varepsilon\rightarrow 0}\limsup\limits_{T\rightarrow\infty}{\diy\frac{1}{T^{m+1}}}
\ln\mathbf{P}(\xi_T(\cdot)\in U_\varepsilon(f))\geq\lim\limits_{\varepsilon\rightarrow 0}\liminf\limits_{T\rightarrow\infty}{\diy\frac{1}{T^{m+1}}}\ln\mathbf{P}(\xi_T(\cdot)\in U_\varepsilon(f))\\
\quad\geq-Q_m\int_0^1f^m(s)ds+\lim\limits_{\varepsilon\rightarrow 0}\liminf\limits_{T\rightarrow\infty}
{\diy\frac{1}{T^{m+1}}}
\ln\mathbf{E}(e^{B_T+N_T(\zeta)\ln2};\zeta_T\in U_\varepsilon(f)).
\end{array} \label{22}\eeq

Observe that in the course if deducing the estimates  (\ref{20}), (\ref{21}) and (\ref{22})
the limit  $T\to \infty$ precedes the limit $\varepsilon \to 0$.

Applying Lemmas \ref{l2.3} and \ref{l2.4} to (\ref{20})--(\ref{22})
completes the proof of the LLDP for the functions  $\in F$. \quad $\blacktriangle$

\begin{remark} \label{r3.1} The above argument allows us to extend the assertion of
Theorem {\rm\ref{th2.1}} to the set of functions $f\in{\rm C}[0,1]$ with $f(0)=0$,
$f(t)\geq 0$ for $0<t\leq 1$ and $f(t)=0$ at finitely many points in $[0,1]$.
\end{remark}

\begin{remark} \label{r3.2} For the Yule process (a process of pure birth  with
$l>0$, $P_l>0$ and $\mu(x)\equiv 0$, see e.g.  \cite{Fell}), the rate functional has the form
$$I(f)=P_l\int_0^1f^l(t)dt, \ \ \ f\in F_M,$$
where $F_M$ is a set of not-decreasing continuous functions $f(t)$ on  $[0,1]$ with $f(0)=0$.
\end{remark}

\medskip

{\sc Proof of Lemma \ref{l2.0}.} Let $N_T(\xi)$ be the number of jumps in process $\xi$
in the time interval  $[0,T]$. In the course of the proof we work on  the event that the
trajectoriy of $\xi$ belongs to $X_T$, i.e., that  $N_T<\infty$. This event has a positive probability.

As was mentioned earlier, the statement of the lemma means that for any measurable set
$G\subseteq X_T$ the equality  \reff{3} is valid.
Denote by  $X_T^{(n)}$ a set of functions $u\in X_T$ with $N_T(u)=n,\ n=0,1,\dots$.
Consider one-to-one mapping
\beq\label{mf}
 u\in X_T^{(n)} \mapsto
 (t_1,\ldots,t_n ; \Delta_1,\ldots ,\Delta_n)\in {\mathfrak X}^{(n)}_T=[0,T]^n_<\times \{+1,-1\}^n,\ n=1,2\ldots .
 \eeq
Here $t_1,\ldots,t_n$ is a sequence of jump times for function $u$ in $[0,T]$,
$\Delta_i$ is a size of jump $u(t_i)-u(t_{i-1})$\ (with $\Delta_1=u(t_1)$). Next $[0,T]^{n}_<$
stands for an $n$-dimensional simplex  $\{(t_1,\ldots,t_n)\ : 0<t_1<\ldots t_n\leq T\}$.

The probabilities $\mathbf{P}(\xi(\cdot)\in G)$ and $\mathbf{P}(\zeta(\cdot)\in G)$ are determiined by

a) the respective densities $\mathfrak{f}_{\xi}$ and $\mathfrak{f}_{\zeta}$ relative to
the summation measure
$\sum\limits_{n\geq 1}\prod\limits_{j=1}^n{\rm d}t_j$
on  ${\mathfrak X}_T:=\bigcup_{n\geq 1}{\mathfrak X}^{(n)}_T$ (here $t_0=0$ as $j=1$), and

b) the probabilities
$\mathbf{P}(\xi(t)=0, 0\leq t\leq T)=e^{-\lambda(0)T}$ ,  $\mathbf{P}(\zeta(t)=0, 0\leq t\leq T)=e^{-T}$.

The densities  $\mathfrak{f}_{\xi}$  and $\mathfrak{f}_{\zeta}$  are of form
\beq\label{fxi}
  \mathfrak{f}_{\xi}(t_1,\ldots,t_n;\Delta_1,\ldots,\Delta_n)
 =\left(\prod_{i=1}^n\nu(x_{i-1},x_i)e^{-h(x_{i-1})\tau_i}\right)e^{-h(x_n)(T-t_n)},\eeq
and
\beq\label{fz}
 \mathfrak{f}_{\zeta}(t_1,\ldots,t_n;\Delta_1,\ldots,\Delta_n)=
2^{-n} \left(\prod_{i=1}^ne^{-\tau_i}\right)e^{-(T-t_n)},
\eeq
where $x_0=0,\ x_i=\sum_{j=1}^i\Delta_j$,
$i=1,\ldots,n$.

Each factor $\nu(x_{i-1},x_i)e^{-h(x_{i-1})\tau_i}$ in (\ref{fxi}) gives the probability density $h(x_{i-1})e^{-h(x_{i-1})\tau_i}$ for the  time the process
 $\xi$ spent at state  $x_{i-1}$ multiplied  by the probability $\nu(x_{i-1},x_i)/h(x_{i-1})$
of a jump from $x_{i-1}$ to $x_i$. The factor  $e^{-h(x_n)(T-t_n)}$
is the probability to stay at  $x_n$ until time $T$. A similar meaning is attributed to the factors
$\frac{1}{2}e^{-\tau_i}$ and $e^{-(T-tto stay a_n)}$. The products of terms in (\ref{fxi}) and (\ref{fz})
reflect the Markovian character of both processes.

The Radon--Nikodym derivative ${\mathfrak p}_T=
{\rm d}{\mathbf P}^{(\xi)}_T(\,\cdot\,\cap X_T)\big/{\rm d}{\mathbf P}^{(\zeta )}_T$ in \eqref{8a}
is a ratio $\mathfrak{f}_\xi /\mathfrak{f}_\zeta $
because the mapping $X^{(n)}_T\to{\mathfrak X}^{(n)}_T$ is one-to-one. The
Radon-Nikodym theorem can be applied here as both densities
$\mathfrak{f}_{\xi}$ and $\mathfrak{f}_{\zeta}$ are positive on  ${\mathfrak X}^{(n)}_T$
$\sum\limits_{n\geq 1}\prod\limits_{j=1}^n{\rm d}t_j$
and measure  ${\mathfrak X}_T$  is finite (for formulation and proof of Radon-Nicodim Theorem
see, e.g., \cite{D-Sch}, Chapter III, Section 10, Theorem 2, or \cite{Ru}, Theorem 6.10). \quad
$\blacktriangle$

\medskip

{\sc Proof of Lemma  \ref{l2.3}.} First, we upper-bound the expected value
$\mathbf{E}(e^{B_T+N_T(\zeta)\ln2};\zeta_T\in U_\varepsilon(f))$.

Given $a>1$ we write:
\beq\label{27}
\begin{array}{l}
\mathbf{E}(e^{B_T+N_T(\zeta)\ln2};\zeta_T\in U_\varepsilon(f)):= E_{1}+E_{2},\\
 E_1:=\mathbf{E}(e^{B_T+N_T(\zeta)\ln2};\zeta_T\in U_\varepsilon(f);N_T(\zeta)\leq T^a),\\
  E_2:=\mathbf{E}(e^{B_T+N_T(\zeta)\ln2};\zeta_T\in U_\varepsilon(f);N_T(\zeta)> T^a).
\end{array}\eeq
Let us bound $E_1$ from above. If $\zeta_T\in U_\varepsilon(f)$ and $N_T(\zeta)\leq T^a$ then,
by virtue of (\ref{1}), it follows that for  any $\gamma_1>0$ and $T$ large enough,
$$\begin{array}{l}B_T=\sum\limits_{i=1}^{N_T(\zeta)}\ln(\nu(\zeta(t_{i-1}),\zeta(t_{i}))\\
\quad\quad\leq T^a(\ln(P_lT^l(M+\varepsilon)^l(1+\gamma_1))+\ln(Q_mT^m
(M+\varepsilon)^m(1+\gamma_1))).\end{array}$$
Here, as before, $M=\max\left[\max\limits_{t\in[0,1]}f(t),1\right]$.

Set $k_1=P_lQ_m(M+\varepsilon )^{l+m}(1+\gamma_1)^2$. Then the following inequality
is fulfilled:
\beq
E_1\leq \exp\{(T^a+1)\ln(k_1T^{l+m})\}.\label{E1}\eeq
\smallskip

Next, we establish an upper bound for  $E_2$.

Fix $\varepsilon$ and $\delta\in (0,1)$ until the completion of the argument.
Denote $M_\delta:=\max\limits_{s\in[0,\delta]}f(s)$. Given  $u\in X_T$, set
$$\widetilde{\nu}(u(t_{i-1}),u(t_i))=\begin{cases}
P_l(u(t_{i-1}))^l,&\hbox{if }\; u(t_i)-u(t_{i-1})=1, \ t_i\geq\delta T,\\
Q_m(u(t_{i-1}))^m,&\hbox{if }\; u(t_i)-u(t_{i-1})=-1, \ t_i\geq\delta T,\\
P_l(T(M_\delta+\varepsilon))^l,&\hbox{if }\;u(t_i)-u(t_{i-1})=1, \ t_i<\delta T,\\
Q_m(T(M_\delta+\varepsilon))^m,&\hbox{if }\ u(t_i)-u(t_{i-1})=-1, \ t_i<\delta T.
\end{cases}$$
As earlier, $t_i$ are the times of jumps in $u$.

If $\zeta_T\in U_\varepsilon(f)$ then, by (\ref{1}) and the form of function
$\widetilde{\nu}(u(t_{i-1}),u(t_i))$, for $T$ sufficiently large  and
$t_{i-1}<\delta T$  we have an inequality
\beq
\widetilde{\nu}(\zeta(t_{i-1}),\zeta(t_i))\geq\nu(\zeta(t_{i-1}),\zeta(t_i)). \label{art4}
\eeq
Next, if $\zeta_T\in U_\varepsilon(f)$ and $\varepsilon$ is sufficiently small then
for $s>\delta$ we have $\zeta_T(s)>\min\limits_{s\in[\delta,1]}f(s)-\varepsilon>0$.
Thus, for $t_{i-1}\geq\delta T$ the condition (\ref{1}) implies that for any $\gamma_2\in (0,1)$
and $T$ large enough
\beq
(1-\gamma_2)\leq\frac{\nu(\zeta(t_{i-1}),\zeta(t_i))}{\widetilde{\nu}(\zeta(t_{i-1}),\zeta(t_i))}\leq(1+\gamma_2).        \label{29}
\eeq
Owing to inequalities (\ref{art4}), (\ref{29}) for any  $\gamma_2\in (0,1)$ and $T$ sufficiently large
we have that
$$\begin{array}{l}
\prod\limits_{i=1}^{N_T(\zeta)}\nu(\zeta(t_{i-1}),\zeta(t_{i}))
\mathbf{1}(\zeta_T\in U_\varepsilon(f),N_T(\zeta)> T^a)\\
\qquad\leq (1+\gamma_2)^{N_T(\zeta)}\prod\limits_{i=1}^{N_T(\zeta)}\widetilde{\nu}(\zeta(t_{i-1}),\zeta(t_{i}))
\mathbf{1}(\zeta_T\in U_\varepsilon(f),N_T(\zeta)> T^a).\end{array}$$

Next, set
$$\widetilde{f}_\delta(s)=\begin{cases}
M_\delta,&\hbox{if }\; t\in[0,\delta),\\
f(s),&\text{if }\; t\in[\delta,1].\end{cases}$$

From the form of  $\widetilde{\nu}(\zeta(t_{i-1}),\zeta(t_i))$ it follows that for
$\zeta_T\in U_{\varepsilon}(f)$ one of inequalities holds true, depending upon the
sign of $\zeta (t_i)-\zeta (t_{i-1})$: either
\beq
\widetilde{\nu}(\zeta(t_{i-1}),\zeta(t_i))\leq P_l(T(\widetilde{f}_\delta(t_{i-1}/T)+\varepsilon))^l,
\label{art1}\eeq
or
\beq
\widetilde{\nu}(\zeta(t_{i-1}),\zeta(t_i))\leq Q_m(T(\widetilde{f}_\delta(t_{i-1}/T)+\varepsilon))^m.
\label{art2}\eeq

If  $\zeta_T\in U_\varepsilon(f)$ then, by virtue of (\ref{10}), process
$\zeta_T$ has \ $\diy\frac{N_T(\zeta)+L}{2}$ \ positive jumps and \
$\diy\frac{N_T(\zeta)-L}{2}$ \ negative jumps. Hence, from (\ref{art1}), (\ref{art2}) we obtain
that
\beq\begin{array}{l}\prod\limits_{i=1}^{N_T(\zeta)}\nu(\zeta(t_{i-1}),\zeta(t_{i}))
\mathbf{1}(\zeta_T\in U_\varepsilon(f),N_T(\zeta)> T^a) \\
\quad \leq\begin{cases}
                              (1+\gamma_2)^{N_T(\zeta)}T^{vL/2}P_l^{\frac{N_T(\zeta)+L}{2}}
Q_m^{\frac{N_T(\zeta)-L}{2}}\prod\limits_{i=1}^{N_T(\zeta)}T^{\frac{l+m}{2}}(M+\varepsilon)^v,
&\hbox{if }\ l\neq m,\\
(1+\gamma_2)^{N_T(\zeta)}P_l^{\frac{N_T(\zeta)+L}{2}}
Q_m^{\frac{N_T(\zeta)-L}{2}}\prod\limits_{i=1}^{N_T(\zeta)}
T^l(\widetilde{f}_\delta(t_{i-1}/T)+\varepsilon)^l,&\text{if }\; l=m.\end{cases}\end{array}
\label{new1}\eeq

Set
$$k_2(T):= \min \left(1,(P_l/Q_m)^{(f(1)-\varepsilon)T/2}\right),\
k_3(T):= \max \left(1,(P_l/Q_m)^{(f(1)+\varepsilon)T/2}\right).$$
Then from (\ref{9}) it follows that
\beq
k_2(T)\leq\biggl(\frac{P_l}{Q_m}\biggl)^{L/2}\leq  k_3(T). \label{31}\eeq

In addition, set $k_4=
\left(\frac{\widetilde{f}_\delta(0)+\varepsilon}{\widetilde{f}_\delta(t_{N_T(\zeta)}/T)+\varepsilon}\right)^l$.
Owing to inequalities \reff{new1}, (\ref{31}), 
for $T$ sufficiently large
$$E_2\leq \begin{cases}
k_3(T)T^{\frac{v(M+\varepsilon)}{2}}\mathbf{E}\prod\limits_{i=1}^{N_T(\zeta)}
                              2P_l^{\frac{1}{2}}
Q_m^{\frac{1}{2}}(1+\gamma_2)T^{(l+m)/2}(M+\varepsilon)^v,&\hbox{if }\; l\neq m,\\
k_3(T)k_4\mathbf{E}\prod\limits_{i=1}^{N_T(\zeta)}2P_l^{\frac{1}{2}}
Q_m^{\frac{1}{2}}(1+\gamma_2)
T^l(\widetilde{f}_\delta(t_i/T)+\varepsilon)^l,&\text{if }\; l=m.
\end{cases}$$

Following Remark  \ref{r4.2} from Appendix, we get an exponential bound for  $E_2$:
$$E_2\leq
\begin{cases} k_3(T)e^{-T}T^{\frac{v(M+\varepsilon)}{2}}\exp\left\{2P_l^{\frac{1}{2}}
Q_m^{\frac{1}{2}}(1+\gamma_2)T^{(l+m)/2+1}(M+\varepsilon)^v\right\},&\hbox{if }\; l\neq m,\\
k_3(T)k_4e^{-T}\exp\left\{2P_l^{\frac{1}{2}}
Q_m^{\frac{1}{2}}(1+\gamma_2)T^{l+1}\int_0^1(\widetilde{f}_\delta(s)+\varepsilon)^lds\right\},
&\hbox{if }\ l=m. \end{cases}$$

Then, for $T$ sufficiently large, selecting  $a<\frac{l+m}{2}+1$ we obtain from the bound
(\ref{E1}) that
$$\mathbf{E}(e^{B_T+N_T(\zeta)\ln2};\zeta_T\in U_\varepsilon(f))=
E_1+E_2\leq 2E_2.$$
Finally, by taking into account that the value $\ln \Big(k_3(T)T^{\frac{v(M+\varepsilon )}{2}}\Big)$
is of order $T\ln T$, while $\ln\big(k_3(T)k_4\big)$ is of order $T$, we conclude:
for any  $\gamma_2\in (0,1)$ and $\delta\in(0,1)$, the following bounds hold true:
$$\begin{array}{l}
\lim\limits_{\varepsilon\rightarrow 0} \limsup\limits_{T\rightarrow\infty}{\diy\frac{1}{T^{v+1}}}
\ln\mathbf{E}(e^{B_T+N_T(\zeta)\ln2};\zeta_T\in U_\varepsilon(f))\\
\qquad\qquad\qquad \leq \begin{cases}0,&\hbox{if }\; l\neq m,\\
2\sqrt{P_lQ_m}(1+\gamma_2)\int_0^1\widetilde{f}_\delta^l(s)ds,&\hbox{if }\ l=m.
\end{cases}\end{array}$$
Taking the limit as $\gamma_2 ,\delta\rightarrow 0$ completes the proof. \quad $\blacktriangle$

\medskip

{\sc Proof of Lemma \ref{l2.4}.} Let us now lower-bound the value  ${E}_2$ from (\ref{27}).
As before, we fix a sufficiently small $\varepsilon$ until the end of the argument. Everywhere
below, $[\,\cdot\,]$ stands for the integer part.

Introduce the event $D:=\left\{\max\limits_{1\leq k\leq N_T(\zeta)+1}\tau_k\leq T^{1-\beta}\right\}$,
where $1<\beta<a$ and $\tau_{N_T(\zeta)+1}:=T-t_{N_T(\zeta)}$. Also consider the event
$C_{\varepsilon}:=\left\{\inf\limits_{t\in [t_{[\varepsilon T/4]},T]}\zeta(t)>\varepsilon/16\right\}$,
where $t_{[\varepsilon T/4]}$ is the time of the  $[\varepsilon T/4]$-th jump in
 $\zeta$.

Obviously
$$\begin{array}{l}
{E}_2 \;=\; 2^{N_T(\zeta)}\mathbf{E}\prod\limits_{i=1}^{N_T(\zeta)}\nu(\zeta(t_{i-1}),\zeta(t_{i}))
\mathbf{1}(\zeta_T\in U_\varepsilon(f),N_T(\zeta)> T^a)\\
\qquad \geq \mathbf{E}2^{N_T(\zeta)}\prod\limits_{i=1}^{N_T(\zeta)}\nu(\zeta(t_{i-1}),\zeta(t_{i}))
\mathbf{1}(D,C_{\varepsilon},\zeta_T\in U_\varepsilon^+(f),N_T(\zeta)> T^a),\end{array}$$
where $U_\varepsilon^+(f):=\left\{g: \min\limits_{t\in[0,1]}g(t)\geq 0\right\}\cap U_\varepsilon(f)$.

Let $\delta=\min\{s:\min\limits_{t\in[s,1]}f(t)\geq 2\varepsilon\}$ and
denote $r(\delta):=\min\{i:t_{i}\geq T\delta\}$.

Suppose that $\zeta_T\in U_\varepsilon^+(f)$ and $r(\delta)+1\leq i\leq N_T(\zeta)$.
By condition (\ref{1}),
depending upon the sign of $\zeta(t_i)-\zeta(t_{i-1})$, we  have ether
\beq
\nu(\zeta(t_{i-1}),\zeta(t_i))\geq(1-\gamma_3)P_l(T(f(t_{i-1}/T)-\varepsilon))^l \label{art5}
\eeq
or
\beq
\nu(\zeta(t_{i-1}),\zeta(t_i))\geq(1-\gamma_3)Q_m(T(f(t_{i-1}/T)-\varepsilon))^m. \label{art6}
\eeq

If the event $C_\varepsilon$ has occurred, and $[\varepsilon T/4]\leq i\leq r(\delta)$, then,
owing to condition (\ref{1}), for any $\gamma_3\in (0,1)$  and a sufficiently large $T$
the following inequality holds true:
\beq
\nu(\zeta(t_{i-1}),\zeta(t_i))\geq(1-\gamma_3)(T\varepsilon/16)^w, \label{art7}
\eeq
where $w:=\min(l,m)$.

For $\zeta_T\in U_\varepsilon^+(f)$ and $1 \leq i\leq [\varepsilon T/4]$ we have
\beq
\nu(\zeta(t_{i-1}),\zeta(t_i))\geq k_5:=\min\left[\inf\limits_{x\in \mathbb{Z}_+}\lambda(x),\inf\limits_{x\in \mathbb{N}}\mu(x)\right]. \label{art8}
\eeq

Let us introduce the function
$$\widehat{f}_\varepsilon(s)=\begin{cases}
\diy\frac{\varepsilon}{16\max(1,P_l,Q_m)},&\hbox{if }\ s\in[0,\delta),\\
f(s)-\varepsilon,&\hbox{if }\ s\in[\delta,1]. \end{cases}$$

Using (\ref{10}), (\ref{art5}), (\ref{art6}), (\ref{art7}), (\ref{art8}), we get the bound
$$\begin{array}{l}
E_2\geq k_6(T)\mathbf{E}\Bigg[ P_l^{\frac{N_T(\zeta)+L}{2}}
Q_m^{\frac{N_T(\zeta)-L}{2}}(1-\gamma_3)^{N_T(\zeta)}2^{N_T(\zeta)}\\
\qquad \times\prod\limits_{i=[\varepsilon T/4]+1}^{N_T(\zeta)}(T\widehat{f}(t_{i-1}/T))^w
\mathbf{1}(D,C_\varepsilon,\zeta_T\in U_\varepsilon^+(f),N_T(\zeta)> T^a)\Bigg],
\end{array}$$
where $k_6(T):=\left(\diy\frac{k_5}{\max(P_l,Q_m)}\right)^{[\varepsilon T/4]}$.

From inequalities (\ref{9}), (\ref{31}) we obtain
$$\begin{array}{l}
E_2\geq k_7(T)\mathbf{E}\bigg[P_l^{\frac{N_T(\zeta)}{2}}
Q_m^{\frac{N_T(\zeta)}{2}}(1-\gamma_3)^{N_T(\zeta)}2^{N_T(\zeta)}\\
\qquad \times\prod\limits_{i=1}^{N_T(\zeta)}(T\widehat{f}(t_{i}/T))^w
\mathbf{1}(D,C_\varepsilon,\zeta_T\in U_\varepsilon^+(f),N_T(\zeta)> T^a)\bigg],\end{array}$$
where  $k_7(T):=\diy\frac{k_6(T)k_2(T)}{M^wT^{w[\varepsilon T/4]}}$.

From Lemma \ref{l4.4} of Appendix it follows that for any
$\gamma_4\in (0,1)$ and $T$ sufficiently large  the following holds true
\beq
{E}_2 \geq k_7(T)\sum\limits_{n=[T^a]+1}^{\infty}2^n(1-\gamma_4)^{n}P_l^{\frac{n}{2}}
Q_m^{\frac{n}{2}}\mathbf{E}
\prod\limits_{i=1}^{n}
T^w(\widehat{f}(t_i/T))^w
\mathbf{1}(D,N_T(\zeta)=n).  \label{35}
\eeq
Here  $\gamma_4$ is expressed via $\gamma_3$ and $\theta$ whereas $\theta\in (0,1)$
is introduced in Lemmas  \ref{l4.3}, \ref{l4.4} from Appendix.

To estimate the product from (\ref{35}), we use Lemma \ref{l4.1}.  Taking into account that
$n>T^a$, we get that  for $T$ large enough,
\beq\begin{array}{l}\mathbf{E}\left[\prod\limits_{i=1}^{n}T^w(\widehat{f} (t_i/T))^w
\mathbf{1}(D,N_T(\zeta)=n)\right]
=\mathbf{E}\left[\prod\limits_{i=1}^{n}T^w(\widehat{f}(t_i/T))^w\mathbf{1}(N_T(\zeta)=n)\right]\\
\qquad\qquad\qquad
-\mathbf{E}\left[\prod\limits_{i=1}^{n}T^w(\widehat{f}(t_i/T))^w\mathbf{1}(\overline{D},N_T(\zeta)=n)\right]\\
\qquad \diy\geq
\frac{\left(T^w\int_0^T(\widehat{f}(t_i/T))^w dt\right)^n}{n!}e^{-T}
-2T^\beta\frac{\left(T^w\int_0^T(\widehat{f}(t_i/T))^w dt-T^{w+1}\alpha_{1/T^\beta}\right)^n}{n!}e^{-T}. \end{array}\label{32}\eeq
Here $\alpha_{1/T^\beta}=\diy\frac{1}{2T^\beta}\inf\limits_{s\in[0,1]}\;\left(\widehat{f}(s)\right)^w=
\frac{1}{2T^\beta}\left(\frac{\varepsilon}{16\max(1,P_l,Q_m)}\right)^w$ (compare to equation
(\ref{47}) in Appendix) .


Let us now estimate the last summand in the right side of (\ref{32}).
Denote $k_8:=\sup\limits_{s\in[0,1]}(\widehat{f}(s))^w$. As $a>\beta$, for a sufficiently large
$T$ the following inequalities hold true:
$$\begin{array}{l}
2T^\beta\left(T^w\int_0^T(\widehat{f}({t}/{T}))^w dt-T^{w+1}\alpha_{1/T^\beta}\right)^n\\
\qquad\leq
2T^\beta\left(T^w\int_0^T(\widehat{f}({t}/{T}))^w dt\right)^n
\left(1-\diy\frac{\varepsilon^w}{2k_8T^\beta (16\max(1,P_l,Q_m))^w}\right)^n\\
\qquad \leq 2T^\beta\left(T^w\int_0^T(\widehat{f}({t}/{T}))^w dt\right)^n
\left(1-\diy\frac{\varepsilon^w}{2k_8T^\beta (16\max(1,P_l,Q_m))^w}\right)^{T^a}\\
\qquad  \leq\left(T^w\int_0^T(\widehat{f}({t}/{T}))^w dt\right)^n
\exp\left(\beta\ln(2T)-\diy\frac{\varepsilon^w}{2k_8(16\max(1,P_l,Q_m))^w}T^{a-\beta}\right)\\
\qquad\leq{\diy\frac{1}{2}}\left(T^w\int_0^T\widehat{f}^w({t}/{T}) dt\right)^n.\end{array}$$
Consequently, from (\ref{32}) it follows that
$$\mathbf{E}\left[\prod\limits_{i=1}^{n}T^w(\widehat{f} (t_i/T))^w
\mathbf{1}(D,N_T(\zeta)=n)\right]\geq\frac{1}{2}
\frac{\left(T^w\int_0^T(\widehat{f}({t}/{T}))^w dt\right)^n}{n!}e^{-T}.$$
By virtue of (\ref{35}), for $T$ sufficiently large,
$$
{E}_2 \geq \frac{k_7(T)}{2}\sum\limits_{n=[T^a]+1}^{\infty}2^n(1-\gamma_4)^{n}\sqrt{P_lQ_m}
\frac{\left(T^w\int_0^T(\widehat{f}({t}/{T}))^w dt\right)^n}{n!}e^{-T}.
$$
From this it follows that, selecting $a<w+1$, for $T$ large enough we obtain the inequalities
\beq \label{33} \begin{array}{l}
E_2 \geq{\diy\frac{k_7(T)e^{-T}}{2}}\exp\left(2(1-\gamma_4)\sqrt{P_l Q_m}
T^{w+1}\int_0^1\widehat{f}(s)^w ds\right)\\
\qquad\qquad
-{\diy\frac{k_7(T)e^{-T}}{2}}\exp\left(a\ln(T)+(w+2)T^a\ln(T)\right)\\
\geq{\diy\frac{k_7(T)e^{-T}}{4}}\exp\left(2(1-\gamma_4)\sqrt{P_l Q_m}
T^{w+1}\int_0^1\widehat{f}(s)^w ds\right).\end{array}\eeq

By virtue of (\ref{33}) and the fact that  $\ln k_7$ is a quantity of order $T\ln T$,
we now obtain that
$$\begin{array}{l} \liminf\limits_{T\rightarrow\infty}{\diy\frac{1}{T^{v+1}}}
\ln\mathbf{E}(e^{B_T+N_T(\zeta)\ln2};\zeta_T\in U_\varepsilon(f))\\
\qquad\qquad\qquad \geq \begin{cases} 0,&\text{if }\; l\neq m,\\
2(1-\gamma_4)\sqrt{P_l Q_m}\int_0^1\widehat{f}(s)^l ds,&\hbox{if }\ l=m.\end{cases}
\end{array}$$
Furthermore, taking into account the definition of function $\widehat{f}(s)$, we obtain that
$$\begin{array}{l}
\lim\limits_{\varepsilon\rightarrow 0}\liminf\limits_{T\rightarrow\infty}{\diy\frac{1}{T^{v+1}}}
\ln\mathbf{E}(e^{B_T+N_T(\zeta)\ln2};\zeta_T\in U_\varepsilon(f))\\
\qquad\qquad\qquad \geq\begin{cases} 0,&\text{if }\; l\neq m,\\
2(1-\gamma_4)\sqrt{P_lQ_m}\int_0^1f^l(s)ds,&\text{if }\; l=m.\end{cases}
\end{array}$$

Taking the limit as $\delta \to 0$ and $\gamma_4\rightarrow 0$  completes the proof of the
lemma. \quad $\blacktriangle$


\section  {Appendix}

\ \ \ \ \ \ In this section we prove the auxiliary assertions used in earlier arguments.

Let $X_T^{(n)}$ stand for the event that process $\zeta$ has exactly $n$ jumps on the
interval $[0,T]$.

\begin{lemma} \label{l4.1} Let  $g(t)$ be a non-negative bounded Borel function and $n\geq 1$.
Then
\beq
\mathbf{E}\left[\prod\limits_{i=1}^n g(t_i)\mathbf{1}(X_T^{(n)})\right]=
\frac{\left(\int_0^Tg(s)ds\right)^n}{n!}e^{-T}, \label{46}
\eeq
and
\beq
\mathbf{E}\left[\prod\limits_{i=1}^ng(t_i)\mathbf{1}(X_T^{(n)})
\mathbf{1}\big(\max\limits_{1\leq k\leq n+1}\tau_k>T\Delta\big)\right]\leq
\frac{2}{\Delta}\frac{\left(\int_0^Tg(s)ds-T\alpha_\Delta\right)^n}{n!}e^{-T}. \label{47}
\eeq
Here  $\Delta>0$ is a constant and $\alpha_\Delta:=\diy\frac{\Delta}{2} \inf\limits_{t\in[0,T]}\;g(t)$.
Further, $t_1,\dots,t_n$ are jump times on $[0,T]$ in process $\zeta$ and $\tau_{n+1}:=T-t_n$.
\end{lemma}

{\sc Proof.} First, we prove \eqref{46}. To this end, write:
$$\mathbf{E}\left(\prod\limits_{i=1}^ng(t_i)\biggl|X_T^{(n)}\right)=
\mathbf{E}\left(\prod\limits_{i=1}^ng(t_i)\biggl|\eta(T)=n\right),$$
where $\eta$ is a Poisson process with mean $\mathbf{E}\eta(t)=t$.

From \cite{KT}, Theorem 2.3, p. 126, it follows that
$$\begin{array}{l}
\mathbf{E}\left(\prod\limits_{i=1}^ng(t_i)\biggl|\eta(T)=n\right)\\
\quad
={\diy\frac{n!}{T^n}}\int_0^T\left(\int_{s_1}^T...\left(\int_{s_{n-1}}^T\prod_{i=1}^n g(s_i)ds_n\right)...ds_2\right)ds_1=
{\diy\frac{1}{T^n}}\left(\int_0^Tg(s)ds\right)^n.\end{array}$$
Therefore,
$$\mathbf{E}\left[\prod\limits_{i=1}^ng(t_i)\mathbf{1}(X_T^{(n)})\right]=
\frac{1}{T^n}\left(\int_0^Tg(s)ds\right)^n\mathbf{P}(\eta(T)=n)=\frac{1}{n!}\left(\int_0^Tg(s)ds\right)^n e^{-T}.
$$

Next, we turn to the proof of \eqref{47}. Here
$$\begin{array}{l}
\mathbf{E}\left[\prod\limits_{i=1}^ng(t_i)\mathbf{1}(X_T^{(n)})
\mathbf{1}\big(\max\limits_{1\leq k\leq n+1}\tau_k>T\Delta\big)\right]\\
\quad \leq\sum\limits_{r=1}^{[2/\Delta]}\mathbf{E}\left[\prod\limits_{i=1}^ng(t_i)\mathbf{1}(\eta(T)=n)
\mathbf{1}\biggl(\eta\biggl(\frac{rT\Delta}{2}\biggl)-\eta\biggl(\frac{(r-1)T\Delta}{2}\biggl)=0\biggl)\right]:=
\sum\limits_{r=1}^{[2/\Delta]}D_r.\end{array} $$
Using the fact that  $\max\limits_{1\leq k\leq n+1}\tau_k>T\Delta$, we get that
there exists an $r$ with $1\leq r \leq [\frac{2}{\Delta}]$ and with no jumps on interval
$\left[\frac{(r-1)T\Delta}{2},\frac{rT\Delta}{2}\right]$.

Write
$$\begin{array}{l}
D_1=\mathbf{E}\left[\prod\limits_{i=1}^ng(t_i)\mathbf{1}(\eta(T)=n)
\mathbf{1}\biggl(\eta\biggl(\frac{T\Delta}{2}\biggl)=0\biggl)\right]\\
\qquad
=\mathbf{E}\left[\prod\limits_{i=1}^ng(t_i)\mathbf{1}\biggl(\eta(T)-\eta\biggl(\frac{T\Delta}{2}\biggl)=n\biggl)
\mathbf{1}\biggl(\eta\biggl(\frac{T\Delta}{2}\biggl)=0\biggl)\right].\end{array}$$
By using the independence of increments in, and the homogeneity of, the Poisson process and
formula (\ref{46}) we obtain
$$D_1=\frac{\left(\int_{\frac{T\Delta}{2}}^Tg(s)ds\right)^n}{n!}e^{-T(1-\Delta/2)}
\mathbf{P}\biggl(\eta\biggl(\frac{T\Delta}{2}\biggl)=0\biggl)=
\frac{\left(\int_{\frac{T\Delta}{2}}^Tg(s)ds\right)^n}{n!}e^{-T}.
$$
Similarly  for any  $1\leq r \leq [\frac{2}{\Delta}]$ one obtains that
$$
D_r=\frac{\biggl(\int\limits_{[0,T]\setminus B_{r,\Delta}}g(s)ds\biggl)^n}{n!}e^{-T},
$$
where $B_{r,\Delta}=\left[\diy\frac{(r-1)T\Delta}{2},\frac{rT\Delta}{2}\right]$.

In view of the relations
$$\min\limits_{1\leq r \leq \left[\frac{2}{\Delta}\right]}\int\limits_{\left[\frac{(r-1)T\Delta}{2},
\frac{rT\Delta}{2}\right]}g(s)ds
\geq T\frac{\Delta}{2}\inf\limits_{s\in[0,T]}g(s)=T\alpha_\Delta ,$$
we get that
$$
\mathbf{E}\left[\prod\limits_{i=1}^ng(t_i)\mathbf{1}(X_T^{(n)})
\mathbf{1}\big(\max\limits_{1\leq k\leq n+1}\tau_k>T\Delta\big)\right]\leq
\frac{2}{\Delta}\frac{\left(\int_0^Tg(s)ds-T\alpha_\Delta\right)^n}{n!}e^{-T}.\blacktriangle
$$

\begin{remark} \label{r4.2} Lemma {\rm\ref{l4.1}} implies that
$$
\mathbf{E}\prod\limits_{i=1}^{\eta(T)}g(t_i)\mathbf{1}(\eta(T)\geq 1)=
e^{-T}\left(\exp\left\{\int_0^Tg(s)ds\right\}-1\right).
$$
\end{remark}


\begin{lemma} \label{l4.3} Consider a sequence  $b_1, b_2, ..., b_n$, where each $b_i$ equals
$-1$ or $1$. Define by  $c_d$ the number of sequences with  following property:
$$
\left|\sum\limits_{k=1}^r b_k\right|\leq d, \ \forall \ 1\leq r \leq n.
$$
Take $d= [T\Delta]$ and $n= O(T^\beta)$ where $T\rightarrow\infty$ while $\Delta>0$, $\beta>1$
are fixed. Then for any $\theta\in (0,1)$ and  sufficiently large $T$ we hjave the bound
$$c_d\geq (1-\theta)^{n+1} 2^n.$$
\end{lemma}

{\sc Proof.} It is clear that if a sequence $b_{2(p-1)d+1},...,b_{2pd}$, with $1 \leq p \leq\diy\frac{n}{2d}$,
has an equal number of  $1$ and $-1$, and in the sequence $b_{2d[\frac{n}{2d}]+1},...,b_n$ the
difference  between the numbers of  $1$ and  $-1$ in the absolute value is at most 1 then required
property is fulfilled.
The number of  such sequences is not less then $\left(C_{2d}^d\right)^{[\frac{n}{2d}]}$.

Using Stirling's formula gives that
$$
\left(C_{2d}^d\right)^{[\frac{n}{2d}]}\sim \left(\frac{\sqrt{2}(2d)^{2d}}{\sqrt{\pi d} d^{2d}}\right)^{[\frac{n}{2d}]}=
\left(\frac{\sqrt{2} 2^{2d}}{\sqrt{\pi d}}\right)^{[\frac{n}{2d}]}\geq 2^{n-2d}(\pi d)^{-\frac{n}{4d}}.
$$
Thus, owing to the fact that $-2d\ln\,2-\diy\frac{nrd}{4d}=o(n)$, we obtain that,
for any $\theta\in (0,1)$ and $T$ sufficiently large,
$$c_d\geq (1-\theta) 2^n \exp\left(-2d \ln2 - \frac{n\ln\pi d}{4d}\right)
\geq (1-\theta)^{n+1} 2^n.\quad \blacktriangle $$


\begin{lemma} \label{l4.4} Take $\beta>1$ and $n\geq T^\beta$ and let $g(\cdot)$ be a non-negative
bounded Borel function. For any  $\theta>0$ and all $T$ sufficiently large the following estimate holds
true
$$\begin{array}{l}
\mathbf{E}\left[g(t_1,...,t_n)\mathbf{1}\left(\max\limits_{1\leq k\leq n+1}\tau_k\leq T^{1-\beta}\right)
\mathbf{1}(X_T^{(n)})\mathbf{1}(\zeta_T\in U_\varepsilon^+(f))\mathbf{1}(C_\varepsilon)\right]\\
\qquad\qquad\qquad \geq(1-\theta)^{2n}\mathbf{E}\left[g(t_1,...,t_n)
\mathbf{1}\left(\max\limits_{1\leq k\leq n+1}\tau_k\leq T^{1-\beta}\right)
\mathbf{1}(X_T^{(n)})\right].\end{array}$$
Here $U_\varepsilon^+(f):=\left\{g: \min\limits_{t\in[0,1]}g(t)\geq 0\right\}\cap U_\varepsilon(f)$,
$C_\varepsilon:=\left\{\inf\limits_{t\in [t_{[\varepsilon T/4]},T]}\zeta(t)>\varepsilon/16\right\}$
and $t_{[\varepsilon T/4]}$ is the point of the $[\varepsilon T/4]$-th jump
in process  $\zeta$.
\end{lemma}

{\sc Proof} As  $f$  is  uniformly continuous on
$[0,1]$, then for $\delta>0$  sufficiently small we have the inequality
$$\sup\limits_{s,t:|s-t|\leq\delta}|f(s)-f(t)|<\frac{\varepsilon}{4}.$$

Fix $\delta$ with $1/\delta\in{\mathbb N}$ and let  $1 \leq r \leq 1/\delta$.

Denote by $B_{m_r,\delta r}$ the event where process $\zeta$ has exactly $m_r$ jumps on
the interval \\
$[T\delta (r-1),T\delta r]$.

Then we can write
$$\begin{array}{l}
\mathbf{E}\left[g(t_1,...,t_n)\mathbf{1}\left(\max\limits_{1\leq k\leq n+1}\tau_k\leq T^{1-\beta}\right)
\mathbf{1}(X_T^{(n)})\mathbf{1}(\zeta_T\in U_\varepsilon^+(f))\mathbf{1}(C_\varepsilon)\right]\\
\;\; =\sum\limits_{m_1,...,m_{1/\delta}}\mathbf{E}\left[g(t_1,...,t_n)
\mathbf{1}\left(\max\limits_{1\leq k\leq n+1}\tau_k\leq T^{1-\beta}\right)
\mathbf{1}(X_T^{(n)})\mathbf{1}(\zeta_T\in U_\varepsilon^+(f))\mathbf{1}(C_\varepsilon)
\prod\limits_{r=1}^{1/\delta}\mathbf{1}(B_{m_r,\delta r})\right].\end{array}$$
The summation here is over all collections with
$\min\limits_{r}m_r\geq \delta T^\beta$, $\sum\limits_r m_r=n$.

Take a collection $m_1,\dots,m_r$ satisfying the above condition.

Consider a piece of a trajectory of $\zeta$ on the interval $[0,\delta T]$.
Denote by $t_{1,1},...,t_{m_1,1}$ the jump points of $\zeta$ lying in this interval.
Suppose that the jumps satisfy the following conditions:

1) The jumps at times $t_{1,1},...,t_{d_{\varepsilon ,1},1}$ are positive where
$d_{\varepsilon ,1}:=[T\varepsilon /4]$.

2) The jumps at times $t_{d_{\varepsilon ,1}+1,1},... ,t_{m_1,1}$ are such that
for any integer $k\in\big[d_{\varepsilon ,1}+1,m_1]$ we have the inequality
$$\left|\sum\limits_{l=d_{\varepsilon ,1}+1}^k \zeta(t_{l,1})\right|\leq
\left[\frac{T\varepsilon}{8}\right].$$

Then, for $T$ large enough, the trajectory $\zeta_T(t)$ has the following properties.

1) The trajectory is non-negative and lies in an $\varepsilon$-neighborhood of function $f$
for $t\in [0,\delta ]$.

2) $\zeta_T(t)\geq\varepsilon/16$ for $t>t_{d_{\varepsilon,1},1}/T$.

3) $|\zeta_T(\delta)-f(\delta)|\leq 3\varepsilon/8$.

Now consider a piece of a trajectory of $\zeta$ defined on the interval $[\delta T,2\delta T]$
and having the property $|\zeta(\delta T)-Tf(\delta )|\leq 3\varepsilon T/8$.
Denote by $t_{1,2},...,t_{m_2,2}$ the jump points of $\zeta$ in this interval.

Let these jumps satisfy the conditions:

1) At times $t_{1,2},...,t_{|d_{\varepsilon,2}|,2}$ the jumps are positive or negative in accordance
with the sign of the value
$d_{\varepsilon,2}:=[T(\max(\varepsilon/4,f(2\delta))-\zeta_T(\delta))]$.

2) At times $t_{|d_{\varepsilon,2}|+1,2},...,t_{m_2,2}$ are such that for any integer
$k\in [|d_{\varepsilon,2}|+1,m_2]$ the following inequality holds true:
$$
\biggl|\sum\limits_{l=|d_{\varepsilon,2}|+1}^k \zeta(t_{l,2})\biggl|\leq \bigg[\frac{T\varepsilon}{8}\bigg].
$$
Then, again for $T$ large enough, the trajectory $\zeta_T(t)$ has the following properties.

1) The trajectory is non-negative and lies in a $\varepsilon$-neighborhood of  $f$ for
$t\in [\delta,2\delta]$.

2) $\zeta_T(t)\geq\varepsilon/16$ as $t\in [\delta,2\delta]$.

3) $|\zeta_T(2\delta)-f(2\delta)|\leq 3\varepsilon/8$.

Further pieces of the trajectory are dealt with by induction.

Let us count the trajectories whose jumps satisfy the above properties.

As $\max\limits_{1\leq k \leq n+1}\tau_k\leq T^{1-\beta}$, we have that for any $r$
with $1 \leq r \leq 1/\delta$ the interval $[T\delta (r-1),T\delta r]$ contains at least
$[\delta T^\beta]$ jumps of process $\zeta$, where $\beta>1$. Using Lemma  \ref{l4.3} yields
that, when  $T$ is sufficiently large, we will have on $[T\delta (r-1),T\delta r]$
a number of pieces of the trajectory with the aforementioned properties which is not
less than
$$(1-\theta)^{m_r+1-|d_{\varepsilon,r}|}
2^{m_r-|d_{\varepsilon,r}|}>(1-\theta)^{2m_r}2^{m_r}.$$

Consequently, the number of trajectories that fulfill the above properties for all $r$
is not less then
\beq
\prod\limits_r(1-\theta)^{2m_r}2^{m_r}=(1-\theta)^{2n}2^n. \label{a.1}
\eeq

Next, the jump directions in $\zeta$ are mutually independent and do not depend either on
the number of jumps within the interval or on the jump times. Hence, we can use equality
\eqref{a.1} and get that
$$\begin{array}{l}
\sum\limits_{m_1,...,m_{1/\delta}}\mathbf{E}\left[g(t_1,...,t_n)\mathbf{1}(\max\limits_{1\leq k\leq n+1}\tau_k\leq T^{1-\beta})\mathbf{1}(X_T^{(n)})\mathbf{1}(\zeta_T\in U_\varepsilon^+(f))\mathbf{1}(C_\varepsilon)
\prod\limits_{r=1}^{1/\delta}\mathbf{1}(B_{m_r,\delta r})\right]\\
\qquad \geq\sum\limits_{m_1,...,m_{1/\delta}}{\diy\frac{(1-\theta)^{2n}2^n}{2^n}}\mathbf{E}
\left[g(t_1,...,t_n)\mathbf{1}(\max\limits_{1\leq k\leq n+1}\tau_k\leq T^{1-\beta})
\prod\limits_{r=1}^{1/\delta}\mathbf{1}(B_{m_r,\delta r})\right]\\
\qquad =(1-\theta)^{2n}\mathbf{E}\Big[g(t_1,...,t_n)
\mathbf{1}(\max\limits_{1\leq k\leq n+1}\tau_k\leq T^{1-\beta})
\mathbf{1}(X_T^{(n)})\Big].\quad\blacktriangle\end{array}$$
\bigskip

\section*{Acknowledgments}

The authors are grateful to A.A. Borovkov, B.M. Gurevich,
A.M. Mogulsky and E.A. Pechersky for their interest and useful comments.
E.A. Pechersky should be credited for an initial statement of the probem.
We also thank the anonymous referees for critical remarks and suggestions.

AVL thanks FAPESP for the financial support via Grant 2017/20482.
YMS thanks Math Department, Penn State University, for hospitality and support
and StJohn's College, Cambridge, for support.
AAY thanks CNPq and  FAPESP for the financial support via
Grants 301050/2016-3 and 2017/10555-0, respectively.


\bigskip

\begin{thebibliography}{99}

\bibitem{DZ} A. Dembo, O. Zeitouni. {\it Large Deviations Techniques and Applications.}
New York, Springer. 1998.

\bibitem{DS}  J.D. Deuschel, D.W.  Stroock. {\it Large Deviations.} AMS Chelsea Publishing,
Providence, RI, 1989.

\bibitem{Holl} F. den Hollander. {\it Large Deviations.} Fields Institute Monographs, {\bf 14}.
American Mathematical Society, Providence, RI, 2000.

\bibitem{OlVar} E. Olivieri, M.E. Vares. {\it Large Deviations and Metastability.} Cambridge
University Press, Cambridge, 2005.

\bibitem{Puh} A. Puhalskii. {\it Large Deviations and Idempotent Probability.} Chapman \& Hall/CRC.
Boca Raton, 2001.

\bibitem{Var}  S.R. Varadhan. {\it Large Deviations and Applications.}  New York, SIAM,1984.

\bibitem{SS1}  Y. Suhov, I. Stuhl. On principles of large deviation and selected data compression. arXiv:1604.06971v1. 2015; Also: I. Stuhl, Y. Suhov,
Selected data compression: a refinement of Shannon's principle In: {\it Analytical and Computational
Methods in Probability  Theory.} Lecture Notes in Computer Science, {\bf 10684}.
Springer, Berlin, 2017, PP. 309--321.

\bibitem{SS2} M. Kelbert, I. Stuhl, Y. Suhov. Weighted entropy and optimal
portfolios for risk-averse Kelly investments. {\it Aequationes Mathematicae}, {\bf 92} No 1 (2018),
165--200.

\bibitem{MSSZ} A. Mazel, Y. Suhov, I. Stuhl, S. Zohren. Dominance of most tolerant species
in multi-type lattice Widom-Rowlinson models. {\it Journ. Stat. Mech.,}  2014. 8--10.

\bibitem{MPY}  A. Mogulsky, E. Pechersky, A. Yambartsev. Large deviations for excursions of
non-homogeneous Markov processes. {\it Electronic Commun. Probab.}, {\bf 19} (2014), 1--8.

\bibitem{VSB} N. Vvedenskaya, Y. Suhov, V. Belitsky. A non-linear model of trading mechanism
on a financial market. {\it  Markov Processes. Rel. Fields,}  {\bf 19} No 1 (2013), 83--98.
arXiv:1201.4580. 2012.

\bibitem{Fell} W. Feller. {\it An Introduction to Probability Theory and Its Applications,}
Vol. {\bf 2}. Wiley, New York, 1971.

\bibitem{KS}
M.Kelbert, Y. Suhov.  {\it Probability and Statistics by Example,} Vol. {\bf 2}. Cambridge University
Press, Cambridge,  2008.

\bibitem{KT} S. Karlin, H. M. Taylor.  {\it A First Course in Stochastic Processes,}
2nd Edition. New York {\it et al.}, Academic Press, 1975.

\bibitem{Korol} V.S. Korolyuk, N.I. Portenko, A.V. Skorokhod, A.F. Turbin. {\it Handbook of
Probability Theory and Mathematical Statistics.} Nauka, Moscow, 1985.

\bibitem{Dyn}
Dynkin E.B. {\it Markov Processes,} vols {\bf 1, 2}.
Academic Press, New York and Springer, Berlin, 1965.

\bibitem{Ito}
K. Ito. {\it Stochastic Processes,} Vol. {\bf 2}. Nauka, Moscow, 1963 (Russian).
Also see K. Ito. {\it Essentials of Stochastic Processes.} Americal Mathematical Society,
Providence, RI, 2006, and  K. Ito. {\it Stochastic Processes.} Springer, Berlin {\it et al.},
2004.

\bibitem{KMcG} S. Karlin, J. McGregor. The classification of birth and death  processes.
{\it Trans. Amer. Math. Soc.}, {\bf 86} (1957), 366--400.

\bibitem{LR} W. Ledermann, G.E.H. Reuter. Spectral theory for the differential equations
of simple birth and death processes.  {\it Philos. Trans. Roy. Soc. London},
Ser. A, {\bf  246}, (1954), 321--369.

\bibitem{No} J.R. Norris. {\it Markov Chains.} Cambridge. Cambridge University Press,
1997.

\bibitem{St} D. Stroock. {\it An Introduction to Markov Processes}, 2nd Ed.
Springer, Heidelberg {\it et al.},  2014.

\bibitem{Bor-Mog2} A.A. Borovkov,  A.A. Mogulskii. On large deviation principles in metric
spaces. {\it Siberian Math. Journ.,} {\bf 51} No 6 (2010), 989--1003.

\bibitem{Bor-Mog1} A.A. Borovkov, A.A. Mogulskii. Large deviation principles for random
walk trajectories I. {\it Theory Probab. Appl.}, {\bf 56} No 4 (2011), 538--561.

\bibitem{Log1} A.V. Logachov. The local principle of large deviations for solutions of
It\^{o} stochastic equations with quick drift. {\it Journ.  Math. Sci.}, {\bf 218} No 1
(2016), 28--38.

\bibitem{Bor-Mog} A.A. Borovkov, A.A. Mogulski.  Inequalities and principle of large
deviations for the trajectories of processes with independent increments.
{\it Siberian Math. Journ.}, {\bf  54} No 2 (2013), 217--226.

\bibitem{D-Sch} N.Dunford, J.Schwartz. {\it Linear operators. General theory.}
Wiley-Blackwell, Hoboken, NJ, 1988.

\bibitem{Ru}
W. Rudin. {\it Real and Complex Analysis.} McGrow-Hill, 1987.

\end{thebibliography}
\end{document}